\documentclass[12pt]{article}
\usepackage{t1enc}
\usepackage[latin1]{inputenc}
\usepackage[french]{babel}

\usepackage{amsmath}
\usepackage{amssymb}

\usepackage{colortbl}

\usepackage[T1]{fontenc} 
\usepackage{times}
\usepackage{graphicx}

\newcommand{\ds}{\displaystyle}

\title{Faustmann Rotation and population dynamics \\
  in the presence of a risk of destructive events} 

\author{ Patrice Loisel\thanks{ INRA, UMR 729 MISTEA, F-34060
    Montpellier, France} \thanks{ SupAgro, UMR 729 MISTEA, F-34060
    Montpellier, France} }

\date{}

\begin{document}

\maketitle

\noindent {\bf Abstract}

The impact of the presence of risk of destructive event on the
silvicultural practice of a forest stand is investigated.  For that,
we consider a model of population dynamics. This model has allowed us
to make the comparison without and with risk, and highlight the
influence of the presence of risk of destructive event on optimal
thinning and optimal rotation period.


\noindent {\large \bf Introduction}

In terms of forest management, the first question that arises is: what
is the optimal rotation period?  In the case where a calculation
method to predict earnings for various rotation period is available, {
  Faustmann (1849)} proposed a formalism based on the expected
discounted income. Many authors have successively improved or
reformulated it ({ [Ohlin ,1921]}, { [Pearse, 1967]}, { [Clark,
  1976]}).  The risk of destruction has been introduced to forest
stands by { Martell (1980)} and { Routledge (1980)} in discrete
time. Thereafter, { Reed (1984)} has studied the optimal forest
rotation in continuous time with the risk of fire.  { Thorsen and
  Helles (1998)} analysed endogeneous risk. { Buongiorno (2001)}
proposed a generalization of Faustmann approach using Markov Decision
Process Models.  { Peyron and Heshmatol-Vezin (2003)} were interested
in natural risks incurred by forests in discrete time. More recently,
{ Goodnow et al. (2008)} took into account, in the stand management,
that the thinning regime affect the proportion of standing trees
damaged for ice-damaged, { Amacher et al. (2008)} highlighted the
influence of silvicultural practices of a landowner on damage if fire
occurs.

{Moreover, many authors have been studying ({ [Näslund,
    1969]}, { [Schreuder, 1971]}, { [Clark,
    1976]}, { [Kao and Brodie, 1980]}, {
    [Roise, 1986]}, { [Haight et al., 1992]}), the
  determination of optimal thinning and cutting age, under certainty,
  using whole-stand models.  More recently, { Hyytiainen
    and Tahvonen (2003)} studied, among other elements, the influence
  of the rate of interest on the rotation period, { Cao et
    al. (2006)} analysed the effect of initial stand states on optimal
  thinning regime and rotation. { Touza et al. (2008)}
  investigated the ecology impact and interaction with the
  management.}

For the absence of risk of destructive events, all the production
cycles are carried out to the same cutting age. When the risk of
destructive event exists and is taken into account, we assume, as most
authors cited above that the operator systematically decides to
interrupt the current cycle and begins a new cycle.  The first
question that follows is about the impact of presence of risk of
destructive event on silviculture. A second question, which is linked
to the first one, is: what is the consequence of silviculture over the
rotation period? To take into account thinning in the calculation of
the land value we used a model of population dynamics.  Contrary to
{ Xabadia and Goetz (2010)} works based on age-structured
models, we consider a simplified model of average individual type i.e.
a model based on characteristics of medium trees, in order to
facilitate interpretation and analyse more specifically the impact of
a risk of destruction on silviculture. The study generalizes easily to
more complex and realistic models.

In a first part, we determine the land value without or with the
presence of risk. In a second part we discuss the silviculture first
in the reference limit case where individual tree growth is
independent of tree-density (no density dependent growth) and then
where individual tree-growth depends on the density (density dependent
growth). Under the assumption of no density dependent growth and for a
fixed period of rotation, we compare the results in the absence and
presence of risk. Then we deduce the impact of the presence of a risk
to forestry. In a third part, under the assumption of density
dependent growth and using the optimal rotation period, we simulate a
{\it Eucalyptus} stand. This study, by taking into account explicitly
implemented thinning, allows us to deepen the results obtained by
{ Reed (1984)} and clarify the underlying assumptions
  (not taking into account thinning and cleanup costs proportional to
  the damage).  The interest to introduce a model of population
  dynamics and to use the proposed method is justified by the
  possibility of taking into account thinning and clearing costs
  depending on the severity of damage therefore.

\noindent {\large \bf The land value}

In the first part we study the land value without the risk of destructive
event through a model of population dynamics. Then we study the same
land value with the risk of destructive event.

\noindent {\bf Without the risk of destructive event}

We first consider a stand without the presence of risk of destructive
event. The study of this case will allow us to define a benchmark
management of the stand.

For a cutting age, or rotation period $T$ and a {rate of} thinning
{by unit of time} $h(.)$, the land value $W_0$ (up to a constant
$c_1$) is the discounted value of cutting incomes net cost of
replanting:

$$ W_0 = \sum_{i=1}^{+\infty} (\mathcal V(h(.),T)-c_1) e^{- i \delta T} =
{(\mathcal V(h(.),T)-c_1) e^{-\delta T} \over 1- e^{-\delta T}} =
{\mathcal V(h(.),T)-c_1 \over e^{\delta T} -1}$$

where $\mathcal V(h(.),T)$ is the income generated by the cutting at
time $T$ and $c_1$ is the cost of replanting.  $\mathcal V(h(.),T)$ is
by definition the sum of the thinning income on period $[0,T]$ and the
income at final cutting age $T$.

The land value can also be interpreted as, $W_0$ is the instantaneous
value of income from the forest in time $ T $ discounted at the
initial time and is solution of:

$$ W_0 = (W_0+\mathcal V(h(.),T)-c_1) e^{-\delta T} $$

\noindent {\bf A model of population dynamic}

To take into account the thinning in the calculation of the land
value, i.e. to express $ \mathcal V(h(.),T) $, we introduce a model of
population dynamics. The considered model is an average tree model:
the state variables are the number $ n $ of trees per hectare and the
{averaged tree-basal area $s$ measured at breast height ($ 1.30$
  meter from the ground)}. The evolution of these two variables is
governed by the system of ordinary differential equations:

\begin{align}
{dn(t) \over dt} &= -(m(t)+h(t)) n(t) \notag \\
{d s(t) \over dt} &= G(n(t),s(t),t) \notag
\end{align}

where $m(.)$ is the natural mortality and $G(.,.,.)$ is the possibly
density dependent growth function: individual tree-growth depends on
the tree density $n(t)$. 

The {considered population dynamic model} permits to explicitly manage
the trees. The use of more {complex and} realistic models with
different classes of tree basal areas, instead of the model presented,
would not cause any additional methodological problem. In order to
facilitate interpretation and analysis of properties and results, we
chose to assume identical the characteristics of the trees.

{\noindent {\bf Link with stand models}

  Such models can be linked with stand models in the following way,
  let the basal area of the stand $S=n.s$ and assume that $n G(n,s,t)
  = g(n,s) \Gamma(t)$ and that $g(n,s) =g_0(S)$ depends only on
  $S=n.s$. Then the dynamics of the basal area is given by:

$${d S(t) \over dt} = g_0(S(t)) \Gamma(t) - (m(t)+h(t)) S(t)$$

In the particular case of $m=0$ and of the price $p$ proportional to
$s$, by substituting first the basal area to the volume and secondly
$h S$ to $h$ we obtain the Clark's model ({ Clark, 1976}
pp 263-269).

\noindent {\bf The forest income}

Once chosen the model of population dynamics we can express the total
income. Total income $\mathcal V(h(.),T)$ for fixed thinning $h(.)$ is
the sum of the thinning income $H(n(.),s(.),h(.),T)$ on the period
$[0,T]$ and the final income $V_0(n(T),s(T))$:

$$ \mathcal V(h(.),T)= H(n(.),s(.),h(.),T) + V_0(n(T),s(T))$$

The thinning income $H(n(.),s(.),h(.),t)$ on $[0,t]$ actualized to
time $t$ is:

$$H(n(.),s(.),h(.),t) = \int_0^t p(s(u)) h(u) n(u) e^{\delta (t-u)} du$$
 
where $n(.)$ and $s(.)$ are solutions of the dynamic model.

\noindent {\bf Calculation of the land value}

The land value is given by:

$$ W_0 = {\mathcal V(h(.),T)-c_1 \over e^{\delta T} -1}$$ 

with: $\ds \mathcal V(h(.),T)= \int_0^T p(s(u)) h(u) n(u) e^{\delta
  (T-u)} du + V_0(n(T),s(T))$

The maximal value of the land value defined by Faustmann is obtained
by solving the problem:

$$({\mathcal P_{F_0}}): \ \ \ \ \ \ \ 
\max_{h(.),T} {\mathcal V(h(.),T)-c_1 \over e^{\delta T} -1}$$

The maximization of the {Faustmann value $W_0$ taking into
  account} thinning $h(.)$ {for} the cutting age $T$ can be
decomposed in two steps: first we maximize $\mathcal V(h(.),T)$ with
respect to $h(.)$:

$$ (\mathcal S_0) \ \ \ \mathcal V_0(T)=\max_{h(.)} \mathcal V(h(.),T)$$

then we maximize $\ds {\mathcal V_0(T)-c_1 \over e^{\delta T} -1}$ with
respect to $T$ with $\mathcal V_0(T)$ resulting from the first step.

\noindent {\bf In the presence of risk of destructive event}

{Following { Reed (1976)}, we suppose that destructive
  events occur in a Poisson process i.e.  that destructive events
  occur independently of one another, and ramdomly in time.

The distribution of the destructive event time is an exponential with
mean $\ds {1 \over \lambda}$: $\ds F(x) = 1- e^{-\lambda x}$ where
$\lambda$ is the expected number of destructive events per unit
time. No assumption is made on the type of destructive
events. Destructive event can also be a storm, a fire, a disease or an
insect infestation.  We assume that the event is always of the same
type.

We assume that $\theta_t$ is the proportion of damaged trees following
a destructive event and $\theta'$ is the rate depreciation of timber
due to an influx of wood on the market. $\theta$ and $\theta'$ are
random variables positively correlated: if $\theta=0$ (no damage) then
$\theta'=0$ and more $\theta$ is high, more $\theta'$ is also. {
  Depending on the type of risk, $\theta_t$ and thus $\theta'_t$ may
  be conditionally distributed to $n(.), s(.)$, to thinning $h(.)$ or
  to another control parameter.} We define the expectations {
  $\alpha$ of $(1-\theta_t)$ and $\alpha_p$ of $(1-\theta_t)(1-
  \theta'_t)$}, due to the assumptions on the random variables
$\theta_t$ and $\theta'_t$ we deduce $\alpha_p \leq \alpha$. {In
  order to lighten the presentation we omit the $n(.),s(.),h(.)$
  dependancy of $\alpha$ and $\alpha_p$. The modelling is suitable for
  different types of destructive risk: the specificity is reflected in
  the choice of random variables $ \theta, \theta'$ distribution laws
  and therefore of functions $ \alpha $ and $\alpha_p$.}

For a cutting age $T$ and a fixed thinning $h(.)$, the land value
$W_0$ is the actualized value at initial time of the sum of two terms.
The first one is {the expectation (with respect to the time of
  event) of} the sum of the land value and the expectation {(with
  respect to $\theta_t$ and $\theta'_t$)} of the total income
$\mathcal V_1(h(.),t)$ for the period $[0,t]$ minus clearing costs
$\mathcal C(n(t),t)$ in case of a destructive event at time $t$. The
second one is the sum of the land value and the total income $\mathcal
V(h(.),T)$ for the period $[0,T]$ in case of no destructive
event. Then $W_0$ is the solution of the following equation:

\begin{equation}
 W_0 = \int_0^T [W_0+\mathcal V_1(h(.),t) -c_1-\mathcal C(n(t),t)]
 e^{-\delta t} dF(t) + (W_0+\mathcal V(h(.),T)-c_1) e^{-\delta T}
 (1-F(T))
\end{equation}

where $\ds \mathcal V_1(h(.),t) = E(\mathcal
V(h(.),\theta,\theta',t))$ with $\mathcal V(h(.),\theta,\theta',t)$
the total income for $\theta$ and $\theta'$, $\mathcal C(n(t),t) =
E(\mathcal C_n(\theta_t,n(t)))$ {with} $\mathcal
C_n(\theta_t,n(t))$ the clearing costs for $\theta_t$.

\noindent {\bf The forest income}

$\ds \mathcal V(h(.),\theta,\theta',t)$ is the sum of the thinning
income $H(n(.),s(.),h(.),t)$ during $[0,t]$ actualized at time $t$ and
the total income $\mathcal V_F(\theta,\theta',t)$ at time $t$:

  $$ \mathcal V(h(.),\theta,\theta',t) = H(n(.),s(.),h(.),t) +
  \mathcal V_F(\theta,\theta',t)$$

  The final income $\mathcal V_F(\theta,\theta',t)$ is assumed
  proportional to the final income without risk of destructive event
  $V_0(n(T),s(T))$ and is given by: $ \ds \mathcal
  V_F(\theta,\theta',t)= (1-\theta)(1-\theta') V_0(n(t),s(t))$.
  From definition of $\alpha$ and $ \alpha_p$ we deduce the total
  income expectation:

$$ \mathcal V_1(h(.),t) = E( \mathcal V(h(.),\theta,\theta',t)) 
= H(n(.),s(.),h(.),t) + \alpha_p(t)V_0(n(t),s(t))$$

In { Reed (1976)}, the total income generated in case of
risk of destructive event is assumed proportional to the income
generated in the case without risk, which is not longer true if one
wants to take the thinning into account.

\noindent {\bf The clearing costs}

{The clearing costs $\mathcal C_n$ are of two types: the first one
  concerns the damaged trees $\theta_t n(t)$ and the second one
  concerns the survival trees $(1-\theta_t) n(t)$.}
The respective contribution of each type of cost depends on the type
of destructive event. Assuming a linear dependency ({
  Amacher et al., 2008}),
we then deduce the expression of the clearing costs:
              
 $$\mathcal C_n(\theta_t,n(t))= c_2+ c_d \theta_t n(t) + c_s (1-\theta_t) n(t)$$

 From the definition of $\alpha$, we deduce the clearing costs
 expectation:

$$ \mathcal C(n(t),t) = E(\mathcal C_n(\theta_t,n(t))) =
c_2+c_n(t) n(t)$$

where $c_n$ is defined by $c_n(t)=c_d(1-\alpha(t)) + c_s
\alpha(t)$. Thus, $\min(c_d,c_s) \leq c_n(t) \leq \max(c_d,c_s)$.

{Except for the fire risk ({ Amacher et al., 1976}), the
  clearing costs mainly concerns the damaged trees}.

Remark: { Reed (1976)} assumed that clearing costs were
independent of the number of damaged trees.  Using a model of
population dynamics allows us to take into account the clearing costs
proportional to the damage.

{\noindent {\bf The salvageable function $\alpha$}

  The salvageable function $\alpha$ dependency is specific of the type
  of risk of destructive event. Hence:

  - for a storm risk, following { Hanewinkel (2008)},
  $\alpha$ is a decreasing function of the ratio height-diameter $ {H
    \over d}$ (where $d$ is the tree-diameter $ d= \sqrt{{4 s \over
      \pi}}$) and the time $t$: then $\alpha$ can be written as $
  \alpha({H(t) \over \sqrt{s(t)}},t)$. More generally $\alpha$ is a
  function of $n(t)$ and $s(t)$.

  - for ice damage risk, following { Goodnow et
    al. (2008)}, $\alpha$ depends on the planting density $n_0$ and on
  the rate of thinning $h(.)$.
 
  - for a fire risk, following { Amacher et al. (2008)},
  $\alpha$ depends on the planting density $n_0$ and on the level of
  preventive intermediate treatment effort $z(.)$.
  In that case, $z(.)$ is an additional control parameter.}

\noindent {\bf Calculation of the land value}

From $(1)$ and expression of $\mathcal V_1(h(.),t)$ we deduce the land
value:

\noindent {\bf Proposition 1}: {\it In the presence of the risk of destructive
  event the land value is given by:
 
$$ W_0 = {\delta
   +\lambda \over \delta}{\mathcal {\widetilde V}_1(h(.),T)-c_1 \over
   e^{(\delta+\lambda) T} -1} -{\lambda \over \delta} (c_1+c_2)$$

 where: $\ds \mathcal {\widetilde V}_1(h(.),T)= \ds \int_0^T [p(s(t)) h(t)
 n(t)+ \lambda \alpha_p(t) V_0(n(t),s(t)) -
 \lambda c_n(t)n(t)] e^{(\delta+\lambda) (T-t)} dt +
 V_0(n(T),s(T))$.}

{\it Proof}: From $(1)$ and expressions of $\mathcal V(h(.),t)$ and
$\mathcal V_1(h(.),t)$:

$$ W_0 = \int_0^T [W_0+H(n(.),s(.),h(.),t) + \alpha_p(t)
V_0(n(t),s(t))-c_1-c_2-c_n(t) n(t)] e^{-\delta t} dF(t)$$

$$+ (W_0+H(n(.),s(.),h(.),T)+V_0(n(T),s(T))-c_1) e^{-\delta T} (1-F(T))$$

We then deduce the land value:

$$ W_0= {\delta +\lambda \over \delta}
{J_0(n(.),s(.),h(.),T) + V_0(n(T),s(T))-c_1 \over
  e^{(\delta+\lambda)T}-1 } -{\lambda \over \delta} (c_1+c_2)$$

where $J_0(n(.),s(.),h(.),T) $ has the following value:

$ \ds \int_0^T [H(n(.),s(.),h(.),t) + \alpha_p(t) V_0(n(t),s(t)) -
c_n(t)n(t) ] e^{\delta (T-t)} d e^{\lambda(T-t)} +
H(n(.),s(.),h(.),T)$.

{After changing the integration order} in $\ds \int_0^T
H(n(.),s(.),h(.),t) e^{\delta (T-t)} d e^{\lambda(T-t)} $, {$J_0$
  becomes}:

 $J_0(n(.),s(.),h(.),T) = \ds \int_0^T [p(s(t)) h(t) n(t)+ \lambda
 \alpha_p(t) V_0(n(t),s(t)) - \lambda c_n(t)n(t)]
 e^{(\delta+\lambda) (T-t)} dt $

then the result. \hfill $\square$

As in the case without risk, we find that the land value can be
deduced from the income $ \mathcal {\widetilde V}_1(h(.),T)$.
{The following differences occur: first the rate of interest $\delta$
  is replaced by $\delta+\lambda$, secondly $\mathcal V(h(.),T)$ is
  replaced by a modified expression of the income $ \mathcal
  {\widetilde V}_1(h(.),T)$.}  {The second difference} is reflected by
the substitution in the case without risk of term
$H(n(.),s(.),h(.),T)$ by the term $J_0(n(.),s(.),h(.),T)$, or { even
  more precisely, the substitution of $\ds p(s(t)) h(t) n(t) $ by $\ds
  p(s(t)) h(t) n(t) +\lambda (\alpha_p(t) V_0(n(t),s(t)) -
  c_n(t)n(t))$.}

The maximal value of the land value is obtained by solving:

$$ (\mathcal P_{F_1}): \ \ \ \max_{h(.),T} W_0= {\delta +\lambda \over \delta}
{\mathcal {\widetilde V}_1(h(.),T)-c_1 \over e^{(\delta+\lambda)T}-1 }
-{\lambda \over \delta} (c_1+c_2)$$

As in the case without risk, the maximal value of the land value with
respect to the thinning $h(.)$ and the cutting age $T$ can be
decomposed in two steps: first we maximize $ \mathcal {\widetilde
  V}_1(h(.),T)$ with respect to $h(.)$:

$$(\mathcal S_1) \ \ \ \mathcal V_1(T) = \max_{h(.)}
  \mathcal {\widetilde V}_1 (h(.),T)  $$
 
then we maximize $\ds{\delta +\lambda \over \delta}
{\mathcal V_1(T)-c_1 \over e^{(\delta+\lambda)T}-1 } -{\lambda \over
  \delta} (c_1+c_2)$ with respect to $T$ with $\mathcal V_1(T)$
resulting from the first step.

\noindent {\large \bf Silviculture for a fixed cutting age $T$ and a no density
  dependent growth}

We consider first the limiting case of no density dependent growth
which will be used as a reference in the case study of density
dependent growth.

\noindent {\bf Without risk of destructive event}

Let us consider the case where individual growth is not density
dependent. In this case the evolution of the tree-{basal area} $ s
$ does not depend on the tree-number $ n $, then does not depend on
the silviculture.

We study the maximization of the land value with respect to thinning
$h(.)$ for a fixed cutting age $T$:

$$(\mathcal S_0) \ \ \ \max_{h(.)} \mathcal V(h(.),T) =  
\int_0^T p(s(t)) h(t) n(t) e^{\delta(T-t)} dt + V_0(n(T),s(T))$$

with the constraint $0 \leq h(t) \leq \overline h$.

Assume that the final income is given by: $V_0(n,s)=p(s) n$. From the
no density dependence of the individual growth, $p(s(t))$ is
independent of the thinnings and only depends on $t$. Then we define
$R(t)= p(s(t))$. We denote the functions $\pi_0$ and $\Pi_0$: $\ds
\pi_0(t) = R'(t)-(\delta+m(t)) R(t)$ and $\Pi_0(t) = \int_t^{T}
e^{\int_u^T (\delta+m(u')+\overline h)du'} \pi_0(u) du$.  Applying the
maximum Pontryagin Principle to the problem $\mathcal P_0$ (see
Appendix A) we can deduce the proposition:

\noindent {\bf Proposition 2}: {\it Assume $\pi_0$ is decreasing, the individual
  tree-growth is not density dependent. Consider a fixed cutting age
  $T$, then the optimal thinnings are given by:

- if $\pi_0(T) \geq 0$ then $h_* \equiv 0$

- if $\pi_0(T) < 0$ then it exists $0 \leq t_* < T$ such that $h_*(t)
=0$ for $t < t_*$ and $h(t)=\overline h$ pour $t > t_*$. {Moreover} if
$\ds \Pi_0(0) > 0$ then $t_*=0$, else $t_*$ is the unique solution of:

$$\Pi_0(t_*)= \int_{t_*}^{T} e^{\int_u^T (\delta+m(u')+\overline h)du'}
\pi_0(u) du =0$$}
  
Remark: The commutation time $t_*$ is a function of the cutting age $T$.

\noindent {\bf In the presence of risk of destructive event}

We study the maximization of the land value with respect to thinning $h(.)$
for a fixed cutting age $T$:

$$ (\mathcal S_1) \ \ \ \max_{h(.)} \mathcal {\widetilde V}_1 (h(.),T) = 
J_0(n(.),s(.),h(.),T) + V_0(n(T),s(T))$$

with the constraint $0 \leq h \leq \overline h$.

{In this paragraph, we consider a function $\alpha$ depending only
  of time $t$. The chosen function $\alpha$ can be interpreted as the
  expectation with respect to the state variables of a more precise
  function depending on the state variables.

  In case of a storm risk, $\alpha$ is a function of the ratio
  height-diameter and time $t$. For a no density dependent growth, the
  tree-basal area $s$ and hence the diameter is an explicit function
  of time $t$, then $\alpha$ depends only on $t$. In this case, the
  approximate function is exact.

Due to the difference in the criterion without and with risk,} 
we deduce that, for the fixed rotation period $T$, the silviculture
differs and depends explicitly on $\lambda, \alpha$ and $\alpha_p$. So
we will pay attention to the consequence for the silvicultural
practice.

We consider, as in the case without risk, a no density dependent
growth for the trees to facilitate the comparaison.

Let denote the functions $\pi_{\lambda}$ and $\Pi_{\lambda}$: $\ds
\pi_{\lambda}(t) = R'(t)-(\lambda(1-\alpha_p(t))+\delta+m(t)) R(t)
-\lambda c_n(t))$ and $\Pi_{\lambda}(t) = \int_t^{T} e^{
  \int_u^T(\lambda+\delta+m(u')+\overline h)du'} \pi_{\lambda}(u) du$.
Applying the maximum Pontryagin Principle to the problem $\mathcal
P_{\lambda}$ (see appendix A)
we deduce the proposition:

\noindent {\bf Proposition 3}: {\it Assume $\pi_{\lambda}$ is decreasing, the individual
  tree-growth is not density dependent.  Consider a fixed cutting age
  $T$, then the optimal thinnings are given by:

- if $\pi_{\lambda}(T) \geq 0$ then  $h_* \equiv 0$

- if $\pi_{\lambda}(T) < 0$ then it exists $0 \leq t_* < T$ such that
$h_*(t) =0$ for $t < t_*$ and $h(t)=\overline h$ for $t > t_*$. 
Moreover if $\ds \Pi_{\lambda}(0) > 0$ then  $t_*=0$, else $t_*$
is the unique solution of:

$$\Pi_{\lambda}(t_*)= \int_{t_*}^{T}
e^{\int_u^T(\delta+\lambda+m(u')+\overline h)du' } \pi_{\lambda}(u) du
=0$$}
  
\noindent {\bf Comparaison: without and with presence of risk}

{By comparing the without and with risk criteria, we remark that
  the rate of interest $\delta$ is replaced by $\delta+\lambda$ and
  the rate of thinning $h(t)$ is replaced by $\ds h(t) +\lambda
  (\alpha_p(t) -{c_n(t) \over p(s(t))})$ if $V_0$ is defined by
  $V_0(n,s)= p(s) n$.

  Moreover, from the definition of $\pi_{\lambda}$, we deduce that
  $\pi_{\lambda}$ is a decreasing function of $\lambda$. Thus, at
  least, in the vicinity of the cutting age $T$, the greater
  $\lambda$, the more frequently $h(T)$ will be equal to $\overline h$
  in the following sense: for $\lambda_1 < \lambda_2$, with the
  associated thinnings $h_1$ and $h_2$, if $h_1(T) =\overline h$ then
  $h_2(T)=\overline h$.}
 
By comparing the two propositions, without considering the clearing
costs, the natural mortality $m(t)$ in the case without risk is
replaced by the mortality due to events $m(t) +\lambda(1-\alpha_p(t))$
in the presence of risk.  It is equivalent also, from a mathematical
point of view, to replace the fixed discount rate $\delta$ by the
variable discount rate $\delta+\lambda(1-\alpha_p(t))$ in the previous
problem.

By comparing the results of the two propositions we deduce that, for a
fixed rotation period $T$, it is usually best to do thinning at least
at the end of the period in the presence of risk. {This comparison
  confirms the previous property deduced from the definition of
  $\pi_{\lambda}$.}

Comparing the results of the two proposals is permitted if the
rotation periods are identical. If we consider the maximization
problem, with respect to the rotation period, the rotation periods
have no reason to be the same. In that case the comparison is not
permitted and only simulations can allow us to compare the respective
thinning. We will therefore perform simulations.

\noindent {\large \bf Silviculture for a density dependent growth}

We now consider the case where individual growth is density dependent.
If the growth is weakly density dependent, by continuity with the case
of no density dependent growth, the obtained results are still valid
at fixed rotation period $T$.  For a greater density dependence, if
$(G(n,s,t)n)'_n$ {is sufficently small} (see Appendix B), the
optimal thinnings are the same as in the previous case in the vicinity
of $T$. If this is not the case, we cannot obtain analytical results
for the solutions, then simulations are required.

We are interested in a stand of {\it Eucalyptus}. {In the absence
  of more precise information on the structure of the salvageable
  function $ \alpha $ for a {\it Eucalyptus} stand, we restricted our
  analysis to $ \alpha $ constant.}  The function of individual growth
is given by:

$$ G(n,s,t) = {0.7445(1-e^{-0.482 n s}) \over n} {dH(t) \over dt}$$

where $H(t)$ is the high at time $t$: $\ds H(t)= H_0 \ds (1- e^{-{t
    \over H_0}})$ and $H_0=30$ the limited high which depends mainly
on soil fertility ({ Saint André et al., 2002}). The
structure of the growth function is generic and can be used for other
species.

{For the clearing costs, we neglected the clearing costs of the
  second type ($c_n(t)=c_d (1-\alpha(t)$)}.

The weight of the trees of {basal area} $s$ and high $h$ is given by:
$\ds v(s,H,t)= 0.29+(127.8+0.32 t)sH $ in kg ({ Saint
  André et al., 2005}). The price is assumed to depend on the weight:
$p(s,t)= 0.1 v(s,H(t),t)-0.25$.

The determination of the cutting age $T$ is important because of its
impact on silviculture. To better describe the silviculture in the
presence of random risk, it is wiser to look at the effective cutting
age $\mathcal T$ and the effective final tree-{basal area}
$\mathcal S$. Thus we calculate the respective expectations and
variances:

\begin{center}
$\ds E(\mathcal T) = \int_0^T t dF(t)+T (1-F(T))= {F(T) \over \lambda}$
\end{center}

$\ds Var(\mathcal T) = \int_0^T (t-{F(T) \over \lambda})^2
  dF(t)+(T-{F(T) \over \lambda})^2 (1-F(T)) ={2 \over
  \lambda^2}(1-F(T))(F(T)-\lambda T)+ {F^2(T) \over \lambda^2}$

$\ds E(\mathcal S) = \int_0^T s(t) dF(t)+ s(T)(1-F(T))$ and $\ds
Var(\mathcal S) = \int_0^T s^2(t) dF(t)+ s^2(T)(1-F(T))$ can be
derived from the simulations.

\noindent {\large \bf Results and Discussion}

The unit of time for the rotation period is the month.  We suppose:
$m=0.0042$ month$^{-1}$, $\lambda=0.0075$ month$^{-1} $, $\delta =
0.0034$ month$^{-1}, \overline h =0.075$ month$^{-1}$.

Assume first, that in case of destructive event, the destruction of
stand is total ($\alpha = 0, \alpha_p=0$). With the constraint of no
thinning, we found (Table 1) the classical well established result:
risk implies a shortening of the optimum rotation period.  By
considering the optimal thinning, the optimal rotation period is
larger than in the case without risk but the expected effective
rotation period is smaller and remains of the same order of magnitude
as in the case without risk.  This can be explained by the fact that,
in the presence of risk, the optimal rotation period $ T $ is achieved
with a relatively low probability: $1-F(T)=e^{-\lambda T}$.  To
complete the study, the standard deviation of the effective rotation
period was calculated. Its value doesn't vary and is about $ 20$
months.  The expected effective final tree-{basal area} varies
only slightly depending on the scenarios.

If the destruction is only partial i.e. a portion of the stumpage is
salvageable ($\alpha=0.6, \alpha_p=0.4$), without risk and for the two
considered tree-density ($n = 650$ or $1650$ stems/ha), there is no
thinning for the optimal solution (Table 2). In presence of risk, with
the constraint of no thinning, the optimal rotation period is close to
the previous one.  In contrast, if we optimize allowing thinning, the
optimal rotation period is greater and thinnings are to be
done. Similarly to the previous case of total destruction, the
expected effective rotation period remains of the same magnitude as in
the case without risk.

{The initial tree density in the studied range did not influence
  the qualitative behavior of optimal management. Comparing Table 1
  and Table 2 for a density of $650$ stems/ha, we remark that, the
  greater the potential damage, the greater the cutting age will be
  and the earlier the beginning of the preventive thinning.  From
  Table 2, we deduced that the greater the initial tree-density, the
  later the beginning of the preventive thinning will be and the lower
  the difference between the land value in presence of risk without
  and with thinning.  In Table 3, for the commonly used cutting age
  value $T=84$ months, we found similar properties for the different
  optimizations.}

 The presence of risk of partial or total destruction involves earlier
thinning. Because of early thinning the amount of standing trees is
smaller at time $ T $. Thus, because of thinning, the rotation period
$ T $, can be extended and greater than the rotation period without
risk. The earlier thinning provides a kind of self-insurance against
risk. We therefore make endogenous the risk through optimization.
 
{ Hyytiainen and Tahvonen (2003)} showed, in another
context under certainty, that a greater rate of interest may lengthen
the optimal rotation for non optimal initial density.  The density of
$650$ stem/ha in not optimal. We have showed that the risk involves in
particular the substitution of the rate of interest by the rate of
interest plus the expected number of events $\lambda$, the consequence
being an increase of the rate of interest. Thus our result is
consistent with the results of Hyytiainen and Tahvonen.}  By observing
the curves of the land value $ W_0 $ (for a density of $650$ stems/ha
and $\alpha=0.6,\alpha_p=0.4$), depending on the rotation problem $ T
$ with optimal thinning, we find that the land value least varies in
the vicinity of the optimal rotation period with risk (Figure 2) than
without risk (Figure 1). This is another consequence of the fact that,
with risk the optimal rotation period $ T $ is achieved with a
relatively low probability.

We consider the case where the rotation period is determined by other
considerations. We take the commonly used value $T = 84$ months (Table
3). Without and with risk, it is optimal to practice 
thinning.  However, in the presence of risk, optimal thinning starts
earlier.

The results obtained in { Reed (1976)} are valid only
under the following assumptions: the manager does not practice
thinning and the clearing costs in case of destructive event are
fixed. The possibility to consider thinning and clearing costs depends
on dammage severity and therefore justifies the interest to introduce
a model of population dynamics.

\noindent {\large{\bf Conclusion}}

We have studied the management of a stand in the presence of risk of
destructive event. In order to determine the optimal thinning relative
to the Faustmann criterion, we have considered a model of population
dynamics, the choosen model is of average tree type. This model has
allowed us to make the comparison without and with risk and
highlighted the influence of the presence of risk of destructive event
on optimal thinning.

Specifically, the obtained land values, without or with the risk of
destructive event, highlighted differences in the criteria to be
maximized. In the case of no density dependent individual growth, we
have highlighted the impact of the presence of risk on the strategies,
generically regardless of the considered species.

In the case of density dependent growth, the calculations for a stand
of {\it Eucalyptus} have shown that the presence of risk of destruction
event involves earlier thinning and a greater rotation period, for the
optimal strategy.

The obtained results are conditioned by the choice of an individual
tree growth model and by the specification of a weight model and a
price model of trees. Other studies using models adapted for other
species would make the obtained results more generic.

\noindent {\bf Acknowledgements} \\
The author would like to thank Olli Tahvonen for helpful discussions
of the subject matter. The author is grateful to the anonymous
reviewers, whose comments and suggestions helped improve the quality
of this paper.

\noindent {\bf References} \\
\noindent { Amacher et al., 2008} Amacher, G.S., Malik, A.S. and
Haight, R.G., Forest landowner decisions and the value of information
under fire risk. {\it Can. J. For. Res.} 35 (2008), pp 2603-2615.

\noindent { Buongiorno, 2001} Buongiorno, J., Generalization of
Faustmann's Formula for Stochastic Forest Growth and Prices with
Markov Decision Process Models. {\it Forest Science}, 47(4) (2001), pp
466-474.

  \noindent { Cao et al., 2006} Cao, T., Hyytiainen, K., Tahvonen,
  O. and Valsta, L., Effects of initial stand states on optimal
  thinning regime and rotation of Picea abies stands. {\it Scandinavian
  Journal of Forest Research}, 21(5) (2006), pp 388-398.

  \noindent { Clark, 1976} Clark, C.W., Mathematical Bioeconomics,
  Wiley, New York (1976).

  \noindent { Faustmann, 1849} Faustmann, M., Berechnung des Wertes
  welchen Waldboden sowie noch nicht haubare Holzbestände für die
  Weldwirtschaft besitzen.  {\it Allgemeine Forst-und Jagd-Zeitung}, 25
  (1849), pp 441-445.

  \noindent { Goodnow at al., 2008} Goodnow, R., Sullivan, J. and
  Amacher, G.S., Ice damage and forest stand management. {\it Journal of
  Forest Economics}, 14(4) (2008), pp 268-288.

\noindent { Haight et al., 1992} Haight, R.G., Monserud,
R.A. and Chew, J.D., Optimal Harvesting with Stand Density Targets:
Managing Rocky Mountain Conifer Stands for Multiple Forest Outputs.
 {\it Forest Science}, 38(3) (1992), pp. 554-574.

  \noindent { Hanewinkel, 2008} Hanewinkel, M., Storm Damage
  Modelling in Southwest Germany based on National Forest Inventory
  Data. The final seminar of the Stormrisk project. Boras, Sweden,
  2008-10-02.

  \noindent { Hyytiainen and Tahvonen, 2003} Hyytiainen, K. and
  Tahvonen, O., Maximum Sustained Yield, Forest Rent or Faustmann:
  Does it Really Matter?  {\it Scandinavian Journal of Forest Research},
  18(5) (2003), pp 457-469.

\noindent { KaoBrodie} Kao, C. and Brodie, J.D., Simultaneous
  Optimization of Thinnings and Rotation with Continuous Stocking and
  Entry Intervals. {\it Forest Science}, Monograph 22 (Supplement to Number
  3) (1980), pp. 338-346.

  \noindent { Martell, 1980} Martell, D.L., The optimal
  rotation of a flammable forest stand, {\it Canadian Journal of Forest
  Research}, 10(1) (1980), pp. 30-34.

\noindent { Näslund, 1969} Näslund, B., Optimal Rotation
and Thinning. {\it Forest Science}, 15(4) (1969), pp. 446-451

  \noindent { Ohlin, 1921} Ohlin, B.,) Concerning the question of
  the rotation period in forestry. {\it Journal of Forest Economics}, vol 1,
  n°1-1995, (1921), pp. 89-114.

\noindent { Pearse, 1967} Pearse, P.H., The optimal forest
rotation.  {\it Forest Chron.}, 43 (1967), pp 178-195.

\noindent { Peyron and Heshmatol-Vazin, 2003} Peyron, J.L. and
Heshmatol-Vazin, M., La modélisation de la forêt landaise:
portée et limites. Communication 7th Conference ARBORA, Pessac, 14-15
Dec 2003, 17p.
 
\noindent { Reed, 1984} Reed, W.J., The Effects of the Risk of Fire on the
Optimal rotation of a Forest. {\it JEEM}, 11 (1984), pp 180-190.

\noindent { Roise, 1986} Roise, J.P., A Nonlinear
Programming Approach to Stand Optimization. {\it Forest Science}, 32(3)
(1986), pp. 735-748.

\noindent { Routledge, 1980} Routledge, R.D., Effect of
potential catastrophic mortality and other unpredictable events on
optimal forest rotation policy. {\it Forest Science}, 26 (1980), pp
386-399.

\noindent { Saint André, 2002} Saint-André, L., Laclau, J.P.,
Bouillet, J.P., Deleporte, P., Mabiala, A., Ognouabi, N., Baillères,
H. and Nouvellon, Y., Integrative modelling approach to assess
the sustainability of the {\it Eucalyptus} plantations in
Congo. Fourth workshop, IUFRO Working Party S5.01.04 « Connection
between Forest Resources and Wood Quality: Modelling Approaches and
Simulation Software». Harrison Hot Springs Resort, British Columbia,
Canada, September 8-15, 2002.

\noindent { Saint André, 2005} Saint-André, L., M'Bou,
A.T., Mabiala, A., Mouvondy, W., Jourdan, C., Roupsard, O., Deleporte,
P., Hamel, O., and Nouvellon, Y., Age-related equations for above- and
below-ground biomass of a {\it Eucalyptus} hybrid in Congo. {\it Forest
Ecology and Management}, 205 (2005), pp 199-214.

\noindent { Schreuder, 1971} Schreuder, G.F., The Simultaneous
  Determination of Optimal Thinning Schedule and Rotation for an
  Even-Aged Forest. {\it Forest Science}, 17(3) (1971), pp. 333-339.

  \noindent { Thorsen and Helles, 1998} Thorsen, B.J. and Helles,
  F., Optimal stand management with endogenous risk of sudden
  destruction. {\it Forest Ecology and Management}, 108(3) (1998), pp. 287-299.

  \noindent { Touza et al., 2008} Touza, J., Termansen, M. and
  Perrings, C., A Bioeconomic Approach to the Faustmann-Hartman Model:
  Ecological Interactions in Managed Forest. {\it Natural Resource
  Modeling}, 21(4) (2008), pp 551-581.

  \noindent { Xabadia and Goetz, 2010} Xabadia, A. and Goetz, R.U.,
  The optimal selective logging time and the Faustmann
  formula. {\it Journal of Forest Economics}, 16 (2010), pp 63-82.

  \noindent {\large \bf Appendix A. Solving the optimisation problem
      $(\mathcal P_{\lambda})$ with a no density dependent growth}

We consider, for a fixed rotation period $T$, the following problem
$(\mathcal P_{\lambda})$:

$$\max_{h(.)} \int_0^T [R(t) h(t) +\lambda (\alpha_p(t) R(t)
- c_n(t))] n(t) e^{(\delta+\lambda)(T-t)}dt + R(T) n(T)$$

with the constraint $0 \leq h(t) \leq \overline h$.
  
We apply the maximum Pontryagin Principle, the Hamiltonian is:

$$H= [R(t) h + \lambda (\alpha_p(t) R(t) -c_n(t) )] n
e^{(\delta+\lambda)(T-t)} - \mu (m(t)+h) n$$

The first order conditions are:

\begin{align} 
& h^*(t) \mbox{  maximizes } \ds [R(t) e^{(\delta+\lambda)(T-t)}- \mu(t)]h n
\mbox{ with } h \in [0,\overline h] \notag \\
\ds {d \mu(t) \over dt} = - {\partial H \over \partial n} = & - [R(t)
h(t)+\lambda (\alpha_p(t) R(t) - c_n(t))]
e^{(\delta+\lambda)(T-t)} + \mu(t) (m(t)+h(t)) \notag 
\end{align}

with the transversality condition $\mu(T) = R(T)$

We consider the evolution of the function $l(t)$ defined by: $\ds l(t)
= R(t) e^{(\delta+\lambda)(T-t)}- \mu(t)$.

$${d l(t) \over dt} = (R'(t)-(\delta+\lambda(1-\alpha_p(t))+m(t)) R(t)
-\lambda c_n(t)) e^{(\delta+\lambda)(T-t)} + l(t)
(m(t)+h(t))$$

and define the function $\pi$ and $\Pi$ by: $\ds \pi(t) =
R'(t)-(\delta+\lambda (1-\alpha_p(t))+m(t)) R(t) -\lambda
c_n(t)$ and $\Pi(t) = \int_t^{T} e^{\int_u^T
  (\delta+\lambda+m(u')+\overline h)du'} \pi(u) du$,

We then deduce the following Proposition:

\noindent {\bf Proposition A.1}: {\it Assume $\pi$ is decreasing, the individual
  tree-growth is not density dependent. Consider a fixed cutting age
  $T$, then the optimal thinnings are given by:

- if $\pi(T) \geq 0$ then  $h_* \equiv 0$

- if $\pi(T) < 0$ then il existe $0 \leq t_* < T$ such that $h_*(t)
=0$ for $t < t_*$ and $h(t)=\overline h$ for $t > t_*$. Moreover if
$\ds \Pi(0) > 0$ then $t_*=0$, else $t_*$ is the unique solution of:

$$\Pi(t_*)= \int_{t_*}^{T} e^{\int_u^T (\delta+\lambda+m(u')+\overline
  h)du'} \pi(u) du =0$$}

{\it Proof:} From $l(T)=0$ and considering the equation in backward
time $t'=T-t$, we deduce that in the vicinity of $T$, $\ds {d l(t')
  \over dt'}$ and $\pi(t)$ are of opposit sign, then $\ds {d l(t)
  \over dt}$ and $\pi(t)$ have the same sign: hence $l(t) < 0$
(resp. $> 0$) if $\pi(t) \geq 0$ (resp. $< 0$). If $\pi(T_*) \geq 0$,
$\pi(t) > 0$ forall $t$ then $l(t) < 0$ forall $t$ which implies $h(t)
=0$.  If $\pi(T_*) < 0$, from $l(t) = -
e^{-\int_{t_*}^T(m(u)+\overline h)du} \Pi(t)$ we deduce that either
$\Pi(0) \leq 0$, $l(t)<0$ and $l(t)$ cannot change of sign or $\Pi(0)
> 0$ and $l(t)$ changes of sign for a unique time $t_*$.

\noindent {\large \bf Appendix B. Solving the optimisation problem
    $(\mathcal P_{\lambda})$ with a density dependent growth}

We consider, for a fixed rotation period $T$, the following problem
$(\mathcal P_{\lambda})$:

$$\max_{h(.)} \int_0^T [p(s(t)) h(t)+\lambda (\alpha_p(t) p(s(t))
- c_n(t))] n(t) e^{(\delta+\lambda)(T-t))}dt + p(s(T))
n(T)$$

with the constraint $0 \leq h(t) \leq \overline h$.
  
We apply the maximum Pontryagin Principle, the Hamiltonian is:

$$H= [p(s) h +\lambda (\alpha_p(t) p(s) -c_n(t))]n
e^{(\delta+\lambda)(T-t)} - \mu_n (m(t)+h) n + \mu_s G(n,s)$$

The first order conditions are:

\begin{align} 
&  h^*(t) \mbox{ maximizes } \ds [p(s(t)) e^{(\delta+\lambda)(T-t)}- \mu(t)]h n
\mbox{ with } h \in [0,\overline h] \notag \\
\ds {d \mu_n(t) \over dt} =  - {\partial H \over \partial n} = &-
[p(s(t)) h(t)+\lambda (\alpha_p(t) p(s(t)) - c_n(t) )]
e^{(\delta+\lambda)(T-t))} + \mu_n(t) (m(t)+h(t)) \notag \\
& -\mu_s(t) G_n(n(t),s(t),t) \notag \\
\ds {d \mu_s(t) \over dt} = - {\partial H \over \partial s} = &-
(p'(s(t)) h(t)+\lambda \alpha_p(t) p'(s(t)))n(t)
e^{(\delta+\lambda)(T-t)} -\mu_s(t) G_s(n(t),s(t),t) \notag 
\end{align}

with the transversality conditions $\mu_n(T) = p(s(T))$ and
$\mu_s(T)=p'(s(T)) n(T)$.

We consider the evolution of the function $l(t)$ defined by: $\ds l(t)
= p(s(t)) e^{(\delta+\lambda)(T-t))}- \mu_n(t)$.

\begin{align} 
{d l(t) \over dt} = & (p'(s(t))
  G(n(t),s(t))-(\delta+\lambda (1-\alpha_p(t))+m(t)) p(s(t)) -\lambda
  c_n(t)) e^{(\delta+\lambda)(T-t)} \notag \\ 
& +\mu_s(t)G_n(n(t),s(t),t) + l(t)(m(t)+h(t)) \notag
\end{align}

which allow us to give the following result:

\noindent {\bf Proposition A.2}: {\it If $p(s) = C s^a$ and $ a (G(n,s,t)n)'_n
  \leq (\delta+\lambda(1-\alpha_p(T))+m(T)) s(0)$ then the thinning
  $h(t)=\overline h$ is optimal in the vicinity of the cutting age
  $T$.}

{\it Proof:} From $l(T)=0$ we deduce: $\ds {d l \over dt}(T) =
-(\delta+\lambda(1-\alpha_p(T))+m(T)) p(s(T)) -\lambda c_n(T)
+p'(s(T))(G(n,s,t)n)'_n(T)$.  From given conditions, $\ds {d l \over
  dt}(T) < 0$ then $l(t) > 0$ in the vicinity of $T$.

{\bf Table 1}: Optimized land value with respect to $T$.
$\alpha=0$, $\alpha_p=0$.

\begin{tabular}{lcccc}
  \hline
Scenario & Optimal thinnings & Cutting age & Land value & Expected effective  \\ 
         &                  &   &            & cutting age \\
         &(month$^{-1}$)& (month)& (euro)  & (month)  \\ 
  \hline
$650$ stems/ha &&&& \\
\hspace{6mm}   Without risk &&&& \\
\hspace{12mm}    $\max_{h(.),T}$ & $h \equiv 0$ & $58.5$ & $2137.5$& $58.5$\\
\hspace{6mm}  With risk &&&& \\
\hspace{12mm} no thinning   &       & $58.5$ & $658.1$ & $47.3$ \\
\hspace{12mm}   $\max_{T}$, no thinning &  & $54$ & $673.3$ & $44.4$\\ 
\hspace{12mm}   $\max_{h(.),T}$ & $h(t)= \overline h, t \geq 36.5$ & $69.5$ & $810.4$ & $54.2$\\
  \hline
\end{tabular}

\vspace{15mm}

{\bf Table 2}: Optimized land value with respect to $T$.
$\alpha=0.6$, $\alpha_p=0.4$.

\begin{tabular}{lcccc}
  \hline
Scenario & Optimal thinnings & Cutting age & Land value & Expected effective  \\ 
         &                  &   &            & cutting age \\
       &   (month$^{-1}$)& (month)& (euro)  & (month)  \\ 
 \hline
$650$ stems/ha &&&& \\
\hspace{6mm}  Without risk& &&& \\
\hspace{12mm}  $\max_{h(.),T}$  & $h \equiv 0$ & $58.5$ & $2137.5$ & $58.5$\\
\hspace{6mm}       With risk &&&& \\
 \hspace{12mm} no thinning        & & $58.5$ & $1108.9$ & $47.3$ \\
\hspace{12mm} $\max_{T}$, no thinning &  & $57.5$ & $1109.8$& $46.7$\\ 
\hspace{12mm} $\max_{h(.),T}$ & $h(t)= \overline h, t \geq 43.5$ & $65.5$ & $1149.0$&$51.7$\\
&&&& \\
 $1650$ stems/ha &&&& \\
\hspace{6mm}   Without risk&&&& \\
\hspace{12mm}   $\max_{h(.),T}$  & $h \equiv 0$ & $59.5$ & $2497.2$& $59.5$\\
\hspace{6mm}       With risk  &&&& \\
 \hspace{12mm}  no thinning        & & $59.5$ & $1230.7$&$48.0$ \\
\hspace{12mm}  $\max_{T}$, no thinning  &  & $58.5$ & $1232.1$& $47.3$ \\ 
\hspace{12mm} $\max_{h(.),T}$ & $h(t)= \overline h, t \geq 46.5$ & $64.5$ & $1251.1$& $51.1$ \\
 \hline
\end{tabular}

\vspace{15mm}
\pagebreak

{\bf Table 3} : Land value for $T=84$ months.

\begin{tabular}{lcccc}
\hline

Scenario & Optimal thinnings & Cutting age & Land value & Expected effective  \\ 
         &                  &   &            & cutting age \\
      &(month$^{-1}$)& (month) & (euro)&(month)  \\ 
\hline
$650$ stems/ha &&&& \\
\hspace{6mm}   Without risk &&&& \\
\hspace{12mm}  $\max_{h(.)}$  &$h(t)=\overline h, t \geq 60$ &$84$ &$1914.1$&$84$\\  
\hspace{6mm}   With risk &&&& \\
\hspace{12mm}    $h(t)=\overline h, t \geq 60$ &     &$84$&  $1018.4$&$62.3$\\
\hspace{12mm}  $\max_{h(.)}$    &   $h(t), t \geq 43.5$&$84$&  $1099.3$&$62.3$ \\
&&&& \\
$1650$ stems/ha &&&& \\
\hspace{6mm}  Without risk &&&& \\
\hspace{12mm}  $\max_{h(.)}$  &$h(t)=\overline h, t \geq 64.5$ &$84$&  $ 2224.0$&$84$ \\ 
\hspace{6mm}  With risk &&&& \\
\hspace{12mm}   $h(t)= \overline h, t \geq 64.5$ &    &$84$&  $1099.3$&$62.3$ \\
\hspace{12mm}  $\max_{h(.)}$    &   $h(t)= \overline h, t \geq 48.5$&$84$&  $1194.9$ &$62.3$\\
\hline
\end{tabular}

\pagebreak

\includegraphics[scale=1.,height=130mm,width=90mm,angle=-90]{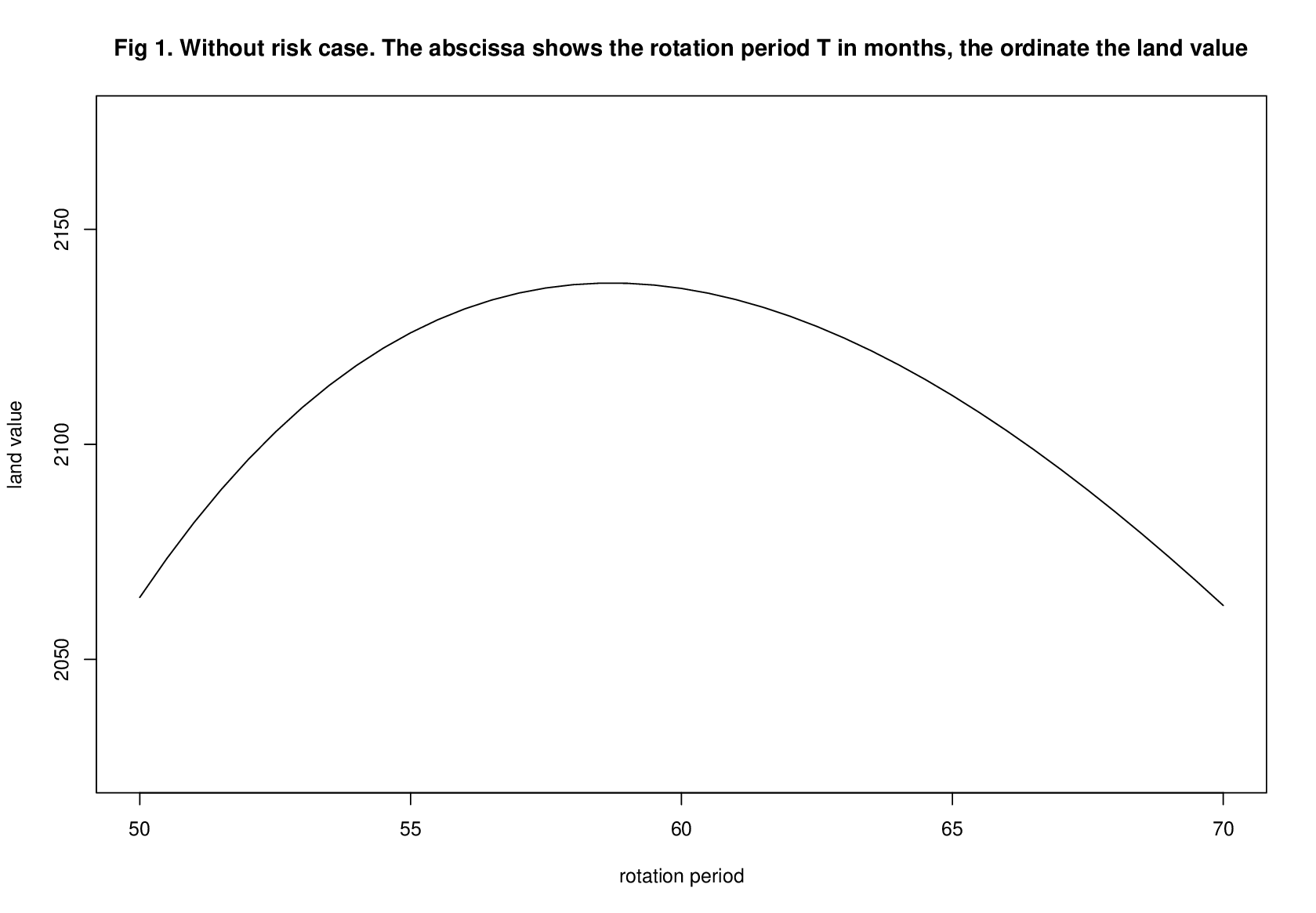}

\includegraphics[scale=1.,height=130mm,width=90mm,angle=-90]{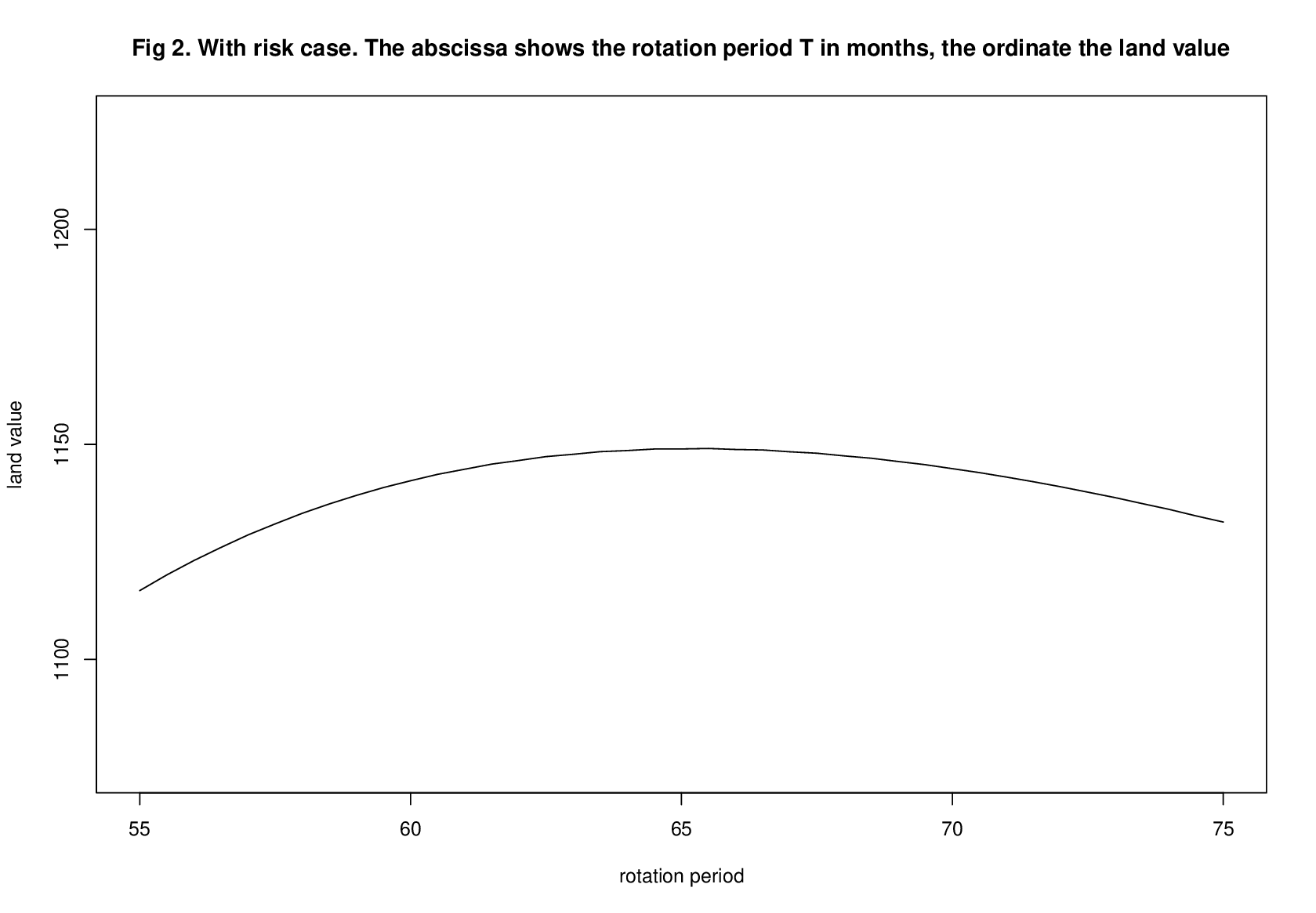}

\end{document}